\documentclass[12pt,a4paper]{article}

\usepackage{amssymb,amsmath}
\usepackage{latexsym}
\usepackage{url}
\usepackage{color}
\usepackage[dvips]{graphicx}

\providecommand{\com[2]}{\begin{tt}[#1: #2]\end{tt}}

\providecommand{\rj[1]}{{\color{red}#1}}
\providecommand{\ad[1]}{{\color{black}#1}}
\newcommand{\re}{{\mathbb R}}
\newcommand{\n}{{\mathbb N}}

\newcommand{\cA}{{\cal{A}}}

\newcommand{\cB}{{\cal{B}}}

\newcommand{\cQ}{{\cal{Q}}}

\newtheorem{theorem}{Theorem}

\newtheorem{lemma}{Lemma}

\newtheorem{remark}{Remark}
\newtheorem{ex}{Example}
\newtheorem{defi}{Definition}
\newtheorem{conj}{Conjecture}

\newtheorem{opq}{Open Question}

\title{Resonance and marginal instability of switching systems}

\author{Vladimir~Y. Protasov and Rapha\"el~M. Jungers 
\thanks{The research of V.P. is supported by the
RFBR grants No
13-01-00642 and 14-01-00332, and by the grant of Dynasty foundation.  R. J. is an F.R.S.-FNRS Research Associate at the ICTEAM Institute of the Universit\'e catholique de Louvain. His work is supported by the Communaut\'e fran\c caise de Belgique - Actions de Recherche
Concert\'ees, and by the Belgian Programme on Interuniversity Attraction Poles initiated
by the Belgian Federal Science Policy Office. }
\thanks{V. P. is with the Department of Mechanics and Mathematics, Moscow State University,
Vorobyovy Gory, 119992, Moscow, Russia. {E-mail: \tt\small
v-protassov@yandex.ru.}}%
\thanks{R. J. is with the ICTEAM Institute,
Universit\'e catholique de Louvain, 4 avenue Georges Lemaitre,
B-1348 Louvain-la-Neuve, Belgium. {E-mail: \tt\small
raphael.jungers@uclouvain.be.}}%
}

\date{}

\begin{document}

\maketitle

\begin{abstract}

We analyse the so-called Marginal Instability of
linear switching systems, both in continuous and discrete time. This is a phenomenon of unboundedness of trajectories
 when the Lyapunov exponent is zero.
 We disprove two recent conjectures of Chitour, Mason, and Sigalotti (2012)
stating that for generic systems, the resonance is sufficient for marginal instability and
for polynomial growth of the trajectories.
The concept of resonance originated with the same authors is modified.
 A characterization of marginal instability under some mild assumptions on the system is provided. These assumptions can be verified algorithmically and are believed to be generic. Finally, we analyze possible types of fastest asymptotic growth of trajectories. An example of a marginally unstable pair of matrices with non-polynomial growth is given.

\smallskip

\smallskip

\begin{flushright}
\noindent  \textbf{AMS 2010} {\em subject
classification: 34F15, 49J15, 39A30}
\end{flushright}

\end{abstract}

\section{Introduction}

We consider the problem of stability for continuous and discrete time linear switching systems
(LSS).
Given a compact family $\cA$ of $d\times d$ matrices, a continuous time LSS is described by the following linear differential
equation on the vector-function $x: \re_+ \to \re^d$:
\begin{equation}\label{cont}
\left\{
\begin{array}{l}
\dot x (t) \ = \ A(t)\, x(t)\,  ;\\
x(0) \, = \, x_0\, ,
\end{array}
\right.
\end{equation}
where $A(\cdot): \, [0, +\infty) \to \cA$ is a measurable function called the {\em switching law}.
The solution $x(\cdot)$ is called the {\em trajectory} of the system corresponding to that switching law and to the
initial condition $x(0) = x_0$. The system is {\em stable} if $x(t) \to 0$ as $t \to \infty$ for every
switching law $A(\cdot)$. Thus, the stability only depends on the family~$\cA$. The measure of stability, called the {\em Lyapunov exponent} $\sigma (\cA),$ is defined as
$$
\sigma (\cA) \ = \ \sup\,  \Bigl( \, \limsup_{t \to \infty} \, \frac1t \,  \log \, \|x(t)\|\, \Bigr) \, ,
$$
where the supremum is computed over all  trajectories of the system. The system is stable if and only if $\sigma (\cA) < 0$.
If $\sigma (\cA) > 0$, then there are unbounded trajectories for which $\|x(t)\|$ grows exponentially as $t \to +\infty$.
In the boundary case $\sigma (\cA) = 0$, the system is never stable, i.e., there is at least one trajectory
that does not converge to the origin as $t \to +\infty$ (see~\cite{MP, B2}). We focus on the question whether the system is {\em bounded} in this case, i.e., all its trajectories are bounded.
\begin{defi}\label{d10}
In the case $\sigma(\cA) = 0$, system~(\ref{cont}) is called marginally stable if it is bounded, otherwise it is called
marginally unstable.
\end{defi}
It is known that the trajectories of marginally unstable systems can grow at most polynomially, and, moreover,
$\|x(t)\| \le C (1+ t^{d-1})\, , \ t \in \re_+$, where $d$ is the dimension of the system. If the family $\cA$ consists of one matrix $A$, then system~(\ref{cont}) becomes
a usual linear ODE $\, \dot x (t)\, = \, A\, x(t)$. Its marginal instability means that the maximal real part of eigenvalues of $A$
is zero and some of them are in resonance, i.e., form a nontrivial Jordan block.

Similarly to~(\ref{cont}), a {\em discrete time switching system} is the following difference equation on the sequence~$\{x(k)\}_{k=0}^{\infty}$:
\begin{equation}\label{dis}
\left\{
\begin{array}{l}
x(k+1) \ = \ A(k)\, x(k)\,  ;\\
x(0) \, = \, x_0\, ,
\end{array}
\right.
\end{equation}
where the switching law $A(k)$ is a sequence of elements from $\cA$. The notions of trajectory, stability, etc., are
extended to discrete systems in a straightforward manner. The stability is measured in terms of the {\em joint spectral radius}
(JSR) of matrices from $\cA$:
$$
\rho(\cA) \ = \ \lim_{k \to \infty} \ \max_{A(j) \in \cA, j = 1, \ldots , k} \ \|A(k)\cdots A(1)\|^{1/k}\,
$$
(in contrast to the Lyapunov exponent, the JSR is usually defined without taking logarithm).
In case of one matrix $\cA = \{A\}$, the JSR becomes the (usual) spectral radius
$\rho(A)$, that is the largest modulus of eigenvalues of~$A$.
The stability of the system is equivalent to the inequality $\rho(\cA) < 1$, and the boundary case $\rho(\cA) = 1$ again
admits the cases of bounded and unbounded systems.  Many tools are available in the literature in order to compute or approximate the joint spectral radius.  See \cite{jsr-toolbox} for a recent toolbox implementing some of them. 
\begin{defi}\label{d20}
In the case $\rho(\cA) = 1$, system~(\ref{dis}) is called marginally stable if it is bounded, otherwise it is called
marginally unstable.
\end{defi}
Marginal stability of continuous and discrete time systems have been studied in the literature
in various contexts. In the theory of refinement functional equations, marginal stability
means that the solutions possess sharp H\"older exponents~\cite{Dub, DL, CH}.
For affine fractal curves, marginal stability is responsible for Lipschitz continuity and
for boundedness of variation~\cite{GS, P3}. It  is also important for  trackability of autonomous agents in sensor networks~\cite{CCJ, J}, in classifications of finite semigroups of integer matrices~\cite{JPB},
in the problem of asymptotic growth  of the Euler binary partition function~\cite{R, P1}
and of some regular sequences~\cite{BCH1, BCH2}. Of course, marginal instability and the
resonance phenomenon are important issues for stability analysis of LSS~\cite{Sun, CMS}.

In the case of one matrix, marginal instability of the discrete system~(\ref{dis}) is equivalent to the fact that $\rho(A) = 1$ and that
some of the largest by modulus eigenvalues of $A$ are in resonance, i.e, form a nontrivial Jordan block.
The fastest growth of trajectories is $\|x(k)\| \asymp k^{\, L}$, where $L$ is the {\em resonance degree}
of the system: the largest size of the block minus one. Some of these conditions are directly extended from one matrix to an arbitrary family~$\cA$. In particular, marginal instability may occur only if the family~$\cA$ is reducible, i.e., in a suitable basis in~$\re^d$, all
matrices of $\cA$ simultaneously get a block upper triangular form (see Definition~\ref{d30} in the next section).
Also, the largest rate of growth (i.e. the degree $L$ of the polynomial) again does not exceed the number of diagonal blocks with the maximal JSR (equal to one)~\cite{P3}. This upper bound was further improved in~\cite{CMS}
by extension of the concept of resonance degree~$L$ to families of matrices.

Necessary conditions for marginal instability were established in~\cite[Theorem 10]{CMS}: the family~$\cA$ has to be reducible and have at least two diagonal blocks in resonance, i.e., they have a common switching law~$A(\cdot)$ which keeps them both away form zero. Thus, similarly to the case of one matrix, {\em the resonance means coincidence of switching laws of several subsystems
providing their fastest growth}. This coincidence generates an unbounded trajectory of the system from bounded ones  of several
subsystems.
However, this is just a necessary condition. Neither reducibility nor resonance implies marginal instability. In fact, finding appropriate sufficient conditions
for marginal instability is a challenging open problem. Two conjectures have been made in~\cite{CMS}. They state that for a generic family~$\cA$ (i.e., for an open and dense subset of matrix families) resonance does imply marginal instability and the fastest rate of growth equals the
resonance degree~$L$.

 {\bf Outline} We start with disproving these two conjectures for both discrete and continuous time systems. This is done in Section~II,
 in Example~\ref{ex10} (discrete time, $d=2$) and in Example~\ref{ex20} (continuous time, $d=4$).
 In these examples, we consider families~$\cA$ of two matrices (so, the convex hull ${\rm co} (\cA)$ is a segment) whose systems are marginally stable, although the resonance degree~$L$ is one, and all admissible perturbations of~$\cA$ respect these properties. This means that a positive resonance degree guarantees neither marginal instability nor polynomial growth
 for generic families. Next, in Section~III, we prove a criterion of marginal instability and a criterion of polynomial
 growth for a wide class of discrete systems (Theorem~\ref{th10}). These are systems such that all irreducible blocks
 of matrices from $\cA$ (in their common block upper-triangular form)
  possess {\em dominant products} (Definition~\ref{d100}). This property can be verified by an algorithm~\cite{GP} and
  is believed to be generic even in a stronger sense: it holds in an open set of matrix families of full Lebesgue measure
  (Conjecture~3). For such systems, the criterion is very simple and easily determines the rate of polynomial growth.
  In particular, the
  resonance phenomenon for such systems means that two blocks have the same dominant word with
  single leading eigenvalues, and those eigenvalues form a Jordan block (Remark~\ref{r5}).

  Finally, in Section~IV, we address the question of possible types of growth for marginally unstable systems. Is it always
  polynomial with an integer degree? This is true for families of one single matrix, for families with dominant products
  (Theorem~\ref{th10}), for families of nonnegative integer matrices~\cite{JPB}. In general, however, the answer is negative as it was shown by Guglielmi and Zennaro~\cite{GZ}. Their examples deal with infinite matrix families~$\cA$, and
  for finite ones the question remained open. We solve it in Theorem~\ref{th20} by  providing examples of pairs of
  $3 \times 3$ matrices for which the fastest growth of trajectories~$x(\cdot)$
  satisfies $\limsup\limits_{k \to \infty} \frac{\log \, \|x(k)\|}{\log k}\, = \, \frac13,$ thus answering open questions from \cite{JPB} (Problems 1 and 2).

\section{Marginal instability of generic systems: two counterexamples}

We start with a short overview of known results on marginal instability. For the sake of simplicity,
we mostly present results for discrete time systems, but unless we mention it explicitly, all the results hold true for continuous time systems as well and
 can be translated in a straightforward way.

 For a compact matrix family $\cA$, we denote
 $$
 M_k \ = \ M_k(\cA) \, = \,
 \max\limits_{A(j) \in \cA, j = 1, \ldots , k} \bigl\| A(k)\ldots  A(1)\bigr\|^{\, 1/k}\, .
 $$
Clearly, $\rho(\cA) = \lim_{k \to \infty} M_k^{1/k}$. It is known~\cite{B1, B2} that if $\rho = 1$, then
$M_k \ge C_1 > 0$, hence the system is not stable (in the sense that there is a trajectory bounded away from the origin). The converse statement requires the notion of irreducibility.
\begin{defi}\label{d30}
A family $\cA$ of $d\times d$-matrices is reducible if all the matrices share a common nontrivial
invariant subspace. Otherwise, the family is said to be irreducible.
\end{defi}
A family is reducible if and only if there is a basis in~$\re^d$ for which all matrices $A \in \cA$
have the following block upper triangular form:
 \begin{equation}\label{blocks}
A \quad = \quad \left(
\begin{array}{cccccc}
A^{(1)} & * &  \ldots &  * \\
0 & A^{(2)}&  * & \vdots \\
\vdots & {} &  \ddots  & * \\
0 &  \ldots &  0 & A^{(r)}
\end{array}
\right)\ ,
\end{equation}
   with irreducible families~$\cA^{(i)} = \{A^{(i)}, \ A \in \cA\}$ in all the diagonal
blocks, $\, i = 1, \ldots , r$. The positions and sizes of all diagonal blocks are the same for all
matrices $A \in \cA$.

If a family $\cA$ with $\rho(\cA)=1$ is irreducible, then it is marginally stable, i.e., there exists a constant $C_2$ such that $\forall  k \in \n\ M_k \le C_2$
(see~\cite{B1, B2}). For reducible families, we always have $M_k \le C_2 k^{d-1}$~\cite{Bell}. This upper bound can be improved
to $C k^{r-1}$, where $r$ is the number of diagonal blocks in factorization~(\ref{blocks}), see~\cite{CH, Sun}.
A further improvement was obtained in~\cite{P2}:  $\, M_k \le C\, k^{r_1-1}$, where $r_1$ is the number of
diagonal blocks~$\cA^{(i)}$ such that $\rho(\cA^{(i)}) = 1$. \\Finally, this estimate was further improved in~\cite{CMS} to
$C k^{L}$, where $L$ is the {\em resonance degree} of the family~$\cA$, which is less by one than the maximal
number of blocks $\rho(\cA^{(i)})$ with the joint spectral radius one that are  {\em in resonance} with each other
(Definition~16 in~\cite{CMS}). Clearly, $L \le r_1-1$ and, to the best of our knowledge, this is the sharpest
estimate known by now. Unfortunately, evaluating $L$ is difficult, because it relies on the ergodic structure of Barabanov norms
of the matrices $\cA^{(i)}$. Let us remember that this norm $\|\cdot\|$ (also called {\em invariant norm}) is characterized by
the property $\max_{A \in \cA^{(i)}} \|Ax\| = \|x\|\, , \ x \in \re^d$. By Barabanov's theorem~\cite{B1} (see also~\cite{B2}
for an analogue for continuous time system), such a norm exists for every irreducible family with joint spectral radius equal to one.

  Hence, each block~$\cA^{(i)}$ has an invariant norm $\|\cdot \|_i$ in the corresponding space $\re^{d_i}.$
The $i$th and $j$th blocks, $i\ne j$, in factorization~(\ref{blocks}) are in {\em resonance}, if there exists a trajectory
$\{x(k)\}_{k \in \{0\}\cup \n}$ such that $\|x_i (k)\|_i = \|x_j (k)\|_j = 1$ for all $k$,
where $x_i$ denotes the restriction of the vector $x$ to the subspace of index $i,$ and $||\cdot||_i$ is a Barabanov norm for the set of blocks of index $i.$
The concept of resonance has been introduced in \cite[definition 6]{CMS}. Thus, resonance formalizes the ``coincidence of switching laws providing
the worst trajectories'' of several subsystems.\\
By theorem~10 of that work, if a system is marginally unstable, then
its family~$\cA$ in factorization~(\ref{blocks}) has at least two blocks in resonance. Then the authors introduced the concept of
resonance degree~$L$ of the system. Roughly, $L$ is the maximal number of blocks simultaneously in resonance (i.e., they are all sharing a common switching law) minus one. The
complete definition is quite technical, and we refer to the original work~\cite{CMS}, definition~16. The resonance degree
gives the best known upper bound for polynomial growth of marginally unstable systems. By theorem~20
of that work, $M_k \, \le \, C \, k^{L}$ for any discrete time system, and a similar estimate  holds for
continuous time systems~\cite[theorem~19]{CMS}. However, these upper bounds are not always attained, and the actual rate of growth of a marginally unstable system may be less than $k^L$. In general, a positive resonance degree does not guarantee
any polynomial growth of trajectories. Nevertheless, it is conjectured in \cite{CMS} that for generic systems, the resonance degree is a correct measure of the polynomial growth. We now formulate these conjectures in slightly different notation.
\begin{conj}\label{con10}~\cite{CMS}
For a generic discrete time system with $\rho = 1$, we have
$$
C_1 k^L \ \le  \ M_k \ \le \ C_2 k^L\, , \quad  k \in \n\, ,
$$
 and the same estimate (replacing $k\in \n\, $ by $\, t \in (0, +\infty)$) holds for generic continuous time systems.
\end{conj}
The second conjecture deals with asymptotic growth of trajectories:
\begin{conj}\label{con20}
Defining $L$ as in Conjecture \ref{con10} above,
for a generic discrete time system with $\rho = 1$, there exists a trajectory $x(k)$
such that $\, \limsup\limits_{k \to \infty} \frac{\log \|x(k)\|}{\log k} = L$,  and
 the same  holds for generic continuous time systems.
\end{conj}

\begin{remark}
{\em These two conjectures are related, and the differences between them deserve an explanation.
The first one is about finite products but claims the property for every single length of products separately, while the second one is only a limit superior, but on the other hand requires the existence of one infinite product with this property, which is stronger.
In both these conjectures the word ``generic'' means that in factorization~(\ref{blocks}) with fixed diagonal  blocks~$\cA^{(i)}$, the set of upper-diagonal entries of all the matrices (which naturally forms a linear space of the corresponding dimension)
contains an open everywhere dense subset such that Conjectures~\ref{con10} and~\ref{con20} hold true, whenever
the set of upper-diagonal elements of matrices from~$\cA$ belong to that subset.}
\end{remark}
\smallskip

We are now presenting two examples that disprove Conjectures~\ref{con10} and~\ref{con20}, for both continuous and discrete switching  systems.

\begin{ex}\label{ex10} (Discrete linear switching systems, $d=2$). {\em We show that there is a pair of two upper-triangular $2\times 2$-matrices~$\cA = \{A_1, A_2\}$, whose joint spectral radius equals one, the diagonal blocks are in resonance, the resonance degree is~$1$,
but the system is marginally stable, and every pair of matrices from an~$\varepsilon$-neighborhood of~$\cA$
(in the two-dimensional space of upper-diagonal entries of matrices $A_1, A_2$) possesses the same property. This shows that the situation, when the rate of polynomial growth of a discrete system equals its resonance degree, is not generic, in contradiction to Conjectures~1 and~2.

Take arbitrary positive numbers~$a$ and~$s$, an arbitrary $2\times 2$-matrix~$B$ such that~$B_{21} = 0$ and all other
  entries are strictly positive, and define
$$
A_1 \ = \
\left(
\begin{array}{cc}
1 & a\\
0& -1
\end{array}
\right)\ ; \qquad
A_2 \ = \ s\, B\, .
$$
Let us show that for any fixed $a\in \re_+$ and~$B\in \re_+^{2\times 2}$ and for all sufficiently small positive~$s$, the pair~$\cA = \{A_1, A_2\}$
possesses all the desired properties. First of all, for the vectors~$u_1 = (1,0)$ and $u_2 = \bigl(1, -\frac{2}{a}\bigr)$,
we have~$A_1u_1 = u_1$ and $A_1u_2 = - u_2$. Whence, the parallelotope~$M$ with vertices~$\pm u_1, \pm u_2$
is invariant with respect to~$A_1$. Consequently,~$A_1$ corresponds to  an isometry in the Minkowski norm~$\|\cdot \|_{M}$ defined by this parallelotope. For every positive~$s < 1/\|B\|_{M}$, we have $\|A_2\|_M < 1$, and therefore,
$\|\cdot \|_M$ is a Barabanov norm for~$\cA$.
 Indeed, for every $x\in \re^2, \, \|x\|_{M}= 1$, we have $\|A_1x\|_M = 1, \|A_2x\|_M < 1$, and so
 $\|x\|_{M} = \max\, \bigl\{ \|A_1x\|\, , \, \|A_2x\|  \bigr\}$. The rate of polynomial growth for~$\cA$ is zero, while the resonance degree is~$1$
(because all infinite products of matrices from~$\cA$ converge to zero, except for those
terminating with  the infinite sequence of~$A_1$). Finally, this property holds after any small enough
perturbation of the elements~$(A_1)_{12} = a$ and $(A_2)_{12}$, because the parallelotope~$M$ and, respectively,
the norm~$\|\cdot \|_{M}$, changes continuously with~$a$.
}
\end{ex}

\begin{ex}\label{ex20} (Continuous linear switching systems, $d=4$). {\em We
present an example of a pair of two~$4\times 4$-matrices, whose continuous time
switching  system possesses the same property
as the discrete time system from  Example~\ref{ex10}: the Lyapunov exponent equals zero, the diagonal blocks are in resonance, the resonance degree is~$1$,
but the system is marginally stable, and so is every pair of matrices from an~$\varepsilon$-neighborhood of~$\cA$.
\smallskip

Take an arbitrary positive $2\times 2$-matrix~$C$,  a positive number~$s$, and an arbitrary block upper triangular $4\times 4$-matrix~$B$,
whose two diagonal $2\times 2$-blocks and all upper-diagonal entries are strictly positive.
Denote
$$
A_1 \ = \
\left(
\begin{array}{cccc}
0 & 0 & C_{11} & C_{12}\\
0& 0 & C_{21} & C_{22}\\
0 & 0 & 0& -1\\
0&0&1&0
\end{array}
\right)\ ; \qquad
A_2 \ = \ B\, - sI.
$$
Let us show that for any fixed $B$ and~$C$, and for all sufficiently large~$s$, the pair~$\cA = \{A_1, A_2\}$
possesses the desired property, i.e., it is marginally stable, although its resonance degree is~$1$.
It is readily seen that the general solution of the system~$\dot x = A_1x$ is
bounded.
 Hence, it possesses a convex  Lyapunov function (a norm)
 $f(x) = \sup\limits_{t \ge 0} \|e^{tA_1}x\|$, which is non-increasing on every trajectory.
 Whence, there is~$\alpha > 0$ such that if~$s$ is large enough, then on every trajectory $x(t)$ of the system~$\dot x = A_2x = (B - sI)x\, , \ x(0) = x_0$ we have~$f(x(t)) < e^{-\alpha t}x(0)$. Thus, the fastest growth of~$f$
 is attained at the constant switching law~$A(t) \equiv A_1$, on which it is bounded. The same is true for
 the system~${\rm co}(\cA)$.
Thus, the rate of polynomial growth for~$\cA$ is zero, while the resonance degree is~$1$.
This property holds after each small
perturbation of the upper right $2\times 2$-blocks of the matrices~$A_1, A_2$. }
\end{ex}

\section{Characterization of marginal instability when blocks have dominant products}
The resonance condition~\cite[definition 6]{CMS} is usually difficult to verify because for that, one needs to know the Barabanov norms of the families~$\cA^{(i)}$, which is notoriously hard to find. Nevertheless, in a vast majority of practical cases, it is possible to decide the marginal instability and to find the rate of polynomial growth. The method is based on the algorithm of exact computation of JSR from~\cite{GP},  on theorem~4 from the same paper~\cite{GP}, and on Theorem~\ref{th10} we are proving  below.

The situation of Examples~\ref{ex10} and~\ref{ex20}, when one matrix of~$\cA$ dominates the others (or, more generally, that a product of these matrices dominates every other product of the same length),
is not that exceptional as it may seem. In fact, it is generally believed that this situation is generic, in the sense that
for almost every matrix family~$\cA$, there is a finite dominant  product of matrices from~$\cA$. To formulate the corresponding conjectures and results
we need to introduce some notation. To every matrix family~$\cA = \{A_1, \ldots , A_m\}$ we associate
an alphabet~$\{a_1, \ldots , a_m\}$, and to an arbitrary word~$\pi$  of this alphabet we assign the corresponding product
$\pi(\cA)$ of matrices
(the empty word corresponds to the identity matrix).
We call a family~$\cA$ {\em normalized} if its JSR is equal to~$1$. To a family~$\cA$ we associate its
normalized family    by dividing each matrix by~$\rho(\cA)$. In what follows we identify a word~$\pi = a_{d_1}\cdots a_{d_n}$ and the corresponding  product~$\pi(\cA) = A_{d_1}\cdots A_{d_n}$. A nonempty word is called simple if is not a power of a shorter word.
\begin{defi}\label{d100} \cite{GP}
A simple word~$\pi$ and the corresponding  product~$\pi(\cA)$ are called dominant for the family~$\cA$ if there is a constant~$q < 1$ such that the spectral radius of each product of
matrices of the normalized family is smaller than~$q$, unless the word of this product is a power of~$\pi$ or of its cyclic permutations.
\end{defi}
The paper~\cite{GP} presents an algorithm for computing the joint spectral radius by ensuring that a given product is dominant. By theorem~4 in~\cite{GP}, the algorithm converges within finite time if and only if the family~$\cA$
possesses a dominant product, whose leading eigenvalue is unique and simple. In case the leading eigenvalue is complex, by uniqueness we mean that there are two conjugate leading eigenvalues. It was shown in~\cite{GP} that
the set of matrices satisfying these assumptions is open and it was conjectured that this set is not only dense
(i.e., generic) but also of the full measure. More precisely, let us denote by $\cQ_m$ the set of all families~$\cA$ of $m$ \, $d\times d$ matrices.
We identify $\cQ_m$ with the space~$\, \re^{\, m\, d^2}$ in the obvious way.
\begin{conj}\label{conj10}
The set of families~$\cA \in \cQ_m$ that possess dominant products with a unique and simple  leading eigenvalue is
of full Lebesgue measure in the set~$\cQ_m$.
\end{conj}
A bulk of numerical examples with matrices from particular applications and
with randomly generated matrices seem to confirm this conjecture (see~\cite{GP}). In all these examples
a dominant product does exist and has a unique and simple leading eigenvalue. Although the corresponding counterexamples are well-known, it is widely believed that they are very rare in practice. Now we are going to see that if each irreducible diagonal
block~$\cA^{(i)}$ possesses this property, then the marginal instability of~$\cA$ can be explicitly decided.
For the sake of simplicity, we formulate the main theorem for the case of two blocks (factorization~(\ref{blocks}) with $r=2$).
\begin{theorem}\label{th10}
Suppose a normalized family of matrices~$\cA$ is reducible and its block upper-triangular form has two blocks~$\cA^{(1)}, \cA^{(2)}$.
Suppose also each block has a dominant word with a simple unique leading eigenvalue; then~$\cA$  is marginally unstable if and only if the following holds:

1) the dominant words of~$\cA^{(1)}$ and of $\cA^{(2)}$ coincide up to   a cyclic permutation;

 2) the leading eigenvalues of the corresponding products $\pi(\cA^{(1)})$ and $\pi(\cA^{(2)})$ are equal, and
the product $\pi(\cA)$ has a nontrivial Jordan block with this eigenvalue.
\smallskip

\noindent  If the conditions 1) and 2) are satisfied, then~$\cA$ has linear growth, i.e., $M_k \asymp k\, , \  k \in \n$.
\end{theorem}
\begin{remark}\label{r5}
{\em Theorem~\ref{th10} suggests a natural concept of resonance for systems with diagonal blocks possessing
dominant words. Two blocks are in resonance if their dominant words coincide up to a cyclic permutation, their leading eigenvalues coincide and form a Jordan block, i.e., one influences the other.
 These conditions are easily verified in practice: first we find dominant words by the algorithm from~\cite{GP},
 compare them and check that the corresponding product~$\pi(\cA)$ has a Jordan block formed by leading eigenvalues
 of the diagonal blocks.}
 \end{remark}
 \begin{remark}\label{r10}
 {\em Conditions 1) and 2) are obviously sufficient for unboundedness, although it remains to prove that the growth of~$\cA$
cannot be faster than linear. Necessity of these conditions, especially of condition 2), is more complicated.
Actually, condition 1) along with the assumption~$\rho(\cA^{(1)}) = \rho(\cA^{(2)})$ means that
the families~$\cA^{(1)}, \cA^{(2)}$
are in resonance in the sense of Definition~6 in~\cite{CMS}. Choosing in an appropriate way
the entries of the upper-right corner~$A_j^{(1, 2)}$ of the matrices~$A_j \in \cA$, we
obtain that the resonance degree of~$\cA$ is one.
However, the
positive resonance degree is still not enough for unboundedness of~$\cA$. The leading eigenvalues of blocks must coincide
(not only their modules, as it holds in the case of resonance~!)  and must form  a Jordan block.}
\end{remark}

We start the proof with several auxiliary results.  Let~$\pi$ be a dominant word for
a family~$\cB$ and $|\pi|=n$. For any natural number~$M,$ we call a word~$x$ {\em white }
if it is represented in the form~$a \pi^k b$, for some words $a$ and $b$ (that may be empty), $|a|\le n-1, |b|\le M$, and $k \in \n \cup \{0\} \cup \{\infty\}$. If $k =\infty$, then $b$ is empty.
 All other words are said to be {\em black}. For example, any word of length less than $n+M$ is white.
 An infinite word is white precisely when it has the form $a \pi^{\infty}$ with $|a| \le n-1$.
\begin{lemma}\label{l10}\cite{GP}
Let a normalized family~$\cB$ possess a dominant word~$\pi$ of length~$n \ge 1$, and the corresponding product
$\pi(\cB)$ have a unique
simple leading eigenvalue. Then there
is $p \in (0,1)$ and $M \ge 0$ such that, for every black word~$\delta$ of length~$\ge n+M$, we have $\|\delta(\cB)\| < p$.
\end{lemma}

\begin{lemma}\label{l20}
For every word $\pi$ of length~$n \ge 1$ and for every number~$M\ge 0$, there is a number~$N \ge M+n$ such that each infinite word
can be split into words so that

a) the first word has length~$\le 2N$;

2) all white words are powers of~$\pi$, and they are separated by black words;

3) the length of each black word, except for the first word (if it is black) is between $n+M$ and $2N$.
\end{lemma}
{\tt Proof.} Let $l$ be the smallest integer such that  $ln > 2n + M$ and let~$N = (l+1)n + M$.
We go along a given  infinite word, starting with the first letter, and find the first power of~$\pi$. We take the maximal power, then go
further starting with the letter following that power, find the first power of~$\pi$, take the maximal power, then go further, etc. This way we spot disjunct powers of~$\pi$. None of those powers can be increased, and there are no extra powers between
the spotted ones. Now write~$\pi^{k_1}, \pi^{k_2}, \ldots , $ for those powers of lengths greater or equal to $nl$, i.e.,
$k_i \ge l$ for all~$i$, and denote by~$x_0$ the word preceding~$\pi^{k_1}$, and by $x_i$ the word between
$\pi^{k_i}$ and $\pi^{k_{i+1}}\, , \, i \ge 1$. Note that~$x_0$ can be empty, but other~$x_i$ cannot.
If there are no powers of~$\pi$ longer than~$nl$, then we have only one infinite word~$x_0$.
If there are finitely many such powers (let~$p^{k_r}$ be the last one), then either~$k_r = \infty$,
or~$x_r$ is an infinite word. Now we make the following transform. For each~$i \ge 0$, we consider the word
$x_i$. If $|x_i| \ge N$, then we split it into words of length between~$N$ and $2N$.
Each of these words is black. Otherwise, if some of them is white, it contains a power of~$\pi$
of length greater than~$N - M - n \ge ln$, i.e., contains the power~$p^{\, l+1}$. Consequently,~$x_i$
contains some of the spotted powers~$\pi^{k_j}$, which is impossible.
Thus, all words $x_i$ longer than~$N$ are split into black words of lengths between~$N$ and~$2N$.
If~$|x_i| < N\, , \, i \ge 1$, then we concatenate~$x_i$ with the word~$\pi$ from the left and with $\pi^{l-1}$ from
the right. We obtain the word~$c = \pi x_i \pi^{l-1}$ which is also black. Indeed, if is is white, then it has
the form~$c = a\pi^kb$, where $|a| \le n-1, |b| \le M$.  The length of~$c$ is greater than~$ln > 2n +M$,
 hence~$k \ge 2$. Therefore, $c$ has prefix~$a\pi^2\, , \, |a| \le n-1$. On the other hand, it has prefix~$\pi x_i$,
 hence it begins with~$\pi^2$, which contradict the construction,because the preceding  word~$\pi^{k_i}$ cannot be extended to the right
  by a power of~$\pi$. Thus, $c$ is black. Its length is smaller than $n + N + (l-1)n < 2N$ and bigger than~$ln +1 > n+M$. This way  we extend every word~$x_i$ of length smaller than~$N$, except for~$x_0$,
to a black word of length between~$n+M$ and~$2N$. After this procedure all the remaining words are powers of~$\pi$, and the lemma follows.

{\hfill $\Box$}
\smallskip

It is known that sets of matrices with a dominant product have a Barabanov norm whose unit ball is a (possibly complex) polytope (see \cite{GP} for more details).  In the proof below, we make use of this norm by default.
This polytope can be obtained by the algorithm from~\cite{GP}.\\

{\tt Proof of Theorem~\ref{th10}}. For each vector $x \in \re^d$, we use the factorization~$x = (x_1, x_2)$,
where the vectors~$x_1, x_2$ correspond to the two blocks. We use the norm~$\|x\| = \|x_1\|_1 + \|x_2\|_2$,
 where here~$\|\cdot \|_i$ is the extremal polytope norm in the subspace of the $i$th block, $i = 1,2$.
 We write~$\pi_1, \pi_2$ for the dominant words of~$\cA^{(1)}, \cA^{(2)}$.

 If the blocks~$\cA^{(1)}$ and~$\cA^{(2)}$ are in resonance, then
there is an infinite word such that, for all finite prefixes of this word, the corresponding products of matrices from~$\cA^{(1)}$
are greater than some constant~$C_0 > 0$, uniformly in the length, and the same is true for the corresponding products of matrices from~$\cA^{(2)}$. Applying Lemma~\ref{l20}
for $\pi = \pi_1$ we split this word. If there are infinitely many black words in that partition, then, by Lemma~\ref{l10},
the norm of products tends to zero, which is impossible. Hence, there are finitely many black products,
and consequently, this word terminates  with~$(\pi_1)^{\infty}$. Applying the same argument for~$\pi = \pi_2$, we see that
this word terminates with~$(\pi_2)^{\infty}$.
Therefore, $\pi_1$ and~$\pi_2$ coincide up to a cyclic permutation. Denote this common dominant word by~$\pi$.
If the corresponding product $\pi$ of matrices from~$\cA$ has no Jordan block with a leading eigenvalue, then~$\|\Pi^r \| \le C , \, r \in \n$, where~$C$ is a constant.\\
We now show that in this case, for every infinite word, all the corresponding partial products (corresponding to the prefixes of this word)
are bounded. Applying Lemma~\ref{l20}, we obtain a partition of this word. Since the length of each
black word is at most~$2N$, it follows that for all products $M_{1,2}$ of matrices from~$\cA$
corresponding to black words, the submatrix $M_{1,2}$ in the upper right block satisfies $\|M_{1,2}x_2 \|_1 \leq C \|x_1\|_1,$
for some common constant~$C$. Take the word $h_k$ that consists of the first $k \ge 5$ words of the partition
and denote $H_k = h_k(\cA)$.
The norms of both diagonal blocks of~$H_k$ are at most~$1$. The rectangular block in the upper right corner
is equal to the following sum:
\begin{equation}\label{eq10}
\sum_{r=1}^kP^{(1)}_1 \cdots P^{(1)}_{r-1} Q_{r}  P^{(2)}_{r+1}\cdots P^{(2)}_k\, ,
\end{equation}
where $P^{(1)}_j$ is the product of matrices from~$\cA^{(1)}$ corresponding to the~$i$th word,
the same is for~$P^{(2)}_j$, and $Q_r$ is the matrix in the upper right corner of the product corresponding to the~$r$th word.
We have $\|Q_r\| \le C$ for all~$r$. Since all black words of the partition, except for the first word (if it is black), are of lengths at least~$n+M$, it follows from Lemma~\ref{l10} that~$\|P^{(i)}_j\|_i < p < 1\, , \, i = 1,2$,
whenever the $j$th word is black and $j\ne 1$.
For all other words, we have
$\|P^{(i)}_j\|_i \le 1\, , \ i =1,2$.   Since there are no two white words in a row, the total number of black words among the first~$k$
words is at least~$\frac{k-1}{2}$. Hence, among all matrices~$P^{(i)}_j$, at least~$\frac{k-3}{2}$
correspond to black words. Therefore, for at least~$\frac{k-5}{2}$ words (we removed the first word),
the norms of the corresponding products are smaller than~$p$, for all other words, they are smaller than or equal to~$1$.
This yields that the norm of each term in~(\ref{eq10}) is bounded above by~$C\, p^{(k-5)/2}$.
So, the total sum is bounded by~$k C\, p^{(k-5)/2}$. Since~$p < 1$ the function~$\varphi(t) = t \, p^{(t-5)/2}$
is bounded above on~$t \in [5, +\infty)$ by some constant~$C_1$. Whence, the norm of the sum~(\ref{eq10}) and, respectively,
the norm of the rectangular block in the upper right corner of~$H_k$ is bounded by~$C_1C$ uniformly in~$k$.

It remains to consider the latter case, when the leading eigenvalues of the two blocks of~$\Pi$ are
equal and form a Jordan block. In this case, this block is of size~$2$, since both these eigenvalues are simple in the
blocks. We have~$\|\Pi^k\| \asymp k$, hence the rate of growth of~$\cA$ is at least linear.
 To show that it is at most linear, we estimate the norm of long products of matrices from~$\cA$ by the same argument as above, with the only
 difference: for white words, we now have~$\|Q_r \| \le C l_r$, where~$l_r$ is the length of the~$r$th word.
 Estimating the same sum~(\ref{eq10}) and summing the lengths, we conclude that the norm of the upper right corner of
 the product~$H_k$ is at most $C\, C_1 |h_k|$, i.e., the norm of~$H_k$ grows not faster than linearly.

{\hfill $\Box$}
\medskip

\section{Sublinear growth}\label{sec-sublin}
What kinds of marginal instability are possible? Is it true that any marginally unstable
discrete system grows polynomially, i.e., $M_N \, \asymp \, N^s$
for some natural number~$s$~? For one matrix $\cA = \{A\}$, the answer is certainly affirmative.
For families of nonnegative integer matrices, it is affirmative as well~\cite{JPB}.
Recently~\cite{BCH2} the affirmative answer was established to arbitrary 
finite families of integer matrices and to some more general families. 
By Theorem~\ref{th10}, the answer is also affirmative for arbitrary families of matrices when all the diagonal blocks $\cA^{(i)}$ have dominant products
with simple leading eigenvalues. This case is believed to be generic (Conjecture~\ref{conj10}).
As we know from the results of~\cite{CMS}, the power $s$ does not exceed
the resonance degree of the system and can be strictly smaller. However, it does not mean that
$s$ is always integer. In fact, the answer is negative, the corresponding example
was constructed by Guglielmi and Zennaro in~\cite{GZ}. This is a family of matrices $\cA = \{A_{\theta}\, | \
\theta \in [0,1]\}$, where
$$
A_{\theta} \ = \
\left(
\begin{array}{cc}
1 & \theta \\
0 & 1 - \theta^p
\end{array}
\right) \, , \qquad \theta \in [0,1]\, .
$$
For any parameter $p \in (1, + \infty)$ this family is marginally unstable and $M_N \asymp N^{p/(1-p)}$.
Thus, even for nonnegative $2\times 2$ matrices, the rate of growth can be sublinear. However, this is an example
of an infinite family, and the question remains about finite ones (see open problems 1 and 2 from~\cite{JPB}). In this section we show  that the sublinear growth
may occur for finite families as well.  Theorem~\ref{th20} provides examples of pairs of $3\times 3$-matrices
$\{A_0, A_1\}$ for which the rate of grows is $N^{1/3}$ and this rate is attained at some trajectory.
This means that, on the one hand, $M_N \le C N^{1/3}$ and, on the other hand, there is a trajectory
$\{x(k)\}_{k =0}^{\infty}$ such that $\, \limsup_{N\to \infty}\frac{\log\, \|x(N)\|}{\log \, N}\, = \, \frac13$.
\smallskip

In this section we use the Euclidean norm by default.
For an arbitrary $\alpha \in \re$, we consider the following pair~$\cA = \{A_0, A_1\}$ of $3\times 3$-matrices
\begin{equation}\label{pr}
A_0 \ = \
\left(
\begin{array}{lll}
1 & 0 & 0\\
0 & 1 & 0\\
0 & 0 & 0
\end{array}
\right)\, , \qquad
A_1 \ = \
\left(
\begin{array}{ccc}
1 & \sin \alpha & \cos \alpha - 1\\
0 & \cos \alpha & -\sin \alpha\\
0 & \sin \alpha &  \cos \alpha
\end{array}
\right)\, .
\end{equation}
\medskip

\begin{theorem}\label{th20}
Let us consider the set $\cA=\{A_0,A_1\}$ described by Equation \eqref{pr} above. If $\, \frac{\alpha }{\pi}\, $ is a quadratic irrational, then there is a constant $C_{\alpha}$ such that
\begin{equation}\label{cubic}
M_N \  \le \  C_{\alpha} \, N^{1/3}\, , \qquad N \in \n \, ,
\end{equation}
 and this rate of growth is attained
 by an infinite product of matrices from~$\cA$.
\end{theorem}

Before we give a proof let us make some observations.  The pair~$\cA$ has two irreducible blocks, of dimensions~$1$ and~$2$.
We denote them by~$\cA^{(1,1)}$ and $\cA^{(2,2)}$, and denote by~$\cA^{(1,2)}$ the $1\times 2$ upper-diagonal block.
Clearly, $\rho(\cA^{(1,1)}) = 1$. Since~$A^{(2,2)}_0$ is an orthogonal projection and~$A^{(2,2)}_1$ is a rotation,
we see that~$\rho(\cA^{(2,2)}) = 1$, and hence~$\rho(\cA) = 1$. The blocks~$\cA^{(1,1)}$ and~$\cA^{(2,2)}$
are in resonance: the sequences~$\{A^n_0\}_{n \in \n}$
and~$\{A^n_1\}_{n \in \n}$ both provide the resonance. It is seen easily that if a product~$\Pi\in \cA^n$
(in what follows, $\cA^n$ denotes the set of products of matrices from~$\cA$ of length~$n$)
has at least two switching points (i.e., cannot be written $A_1^sA_0^{n-s}$ or  $A_0^sA_1^{n-s}$), then $\|\Pi^{(2,2)}\| < 1$. Indeed, the matrix~$A_0^{(2,2)}$
projects the unit ball of~$\re^2$ on the segment~$[-1, 1]$, while any power of the rotation~$A_1^{(2,2)}$ turns this segment
to a nonzero angle (because~$\alpha/\pi$ is irrational).  Thus, after the second switching, the projection~$A_0^{(2,2)}$
  maps this segment strictly inside the unit ball, and the norm of the product is smaller than~$1$.
  As a consequence, not only all resonance products but also all extremal products (whose norms are equal to one identically for all
  suffices) have at most two switches. In the proof below we shall see that all such products
  of matrices from~$\cA$ are actually bounded, and the infinite growth appears only for products with
   an infinite number of switches. Thus, rather surprisingly, for the matrices in (\ref{pr}), {\em the infinite growth is attained, but not on extremal products}.
\smallskip

In the proof of Theorem~\ref{th20} we will need the following technical lemma
\begin{lemma}\label{l30}
For any~$p \ge 2\ $ and~$\ t \in \bigl(0, \frac{\pi}{2}\bigr]$, we have
\begin{equation}\label{sin-cos}
\sin t \, + \, p \, \cos t \ \le \ \sqrt[3]{p^3 + \frac{20}{\sin t}}.
\end{equation}
\end{lemma}
{\tt Proof}. We consider two cases.

$\mathbf{1}.\  \sin t \ge \frac{2}{p}.$ In this case we prove a stronger inequality: $\sin t + p \cos t\, \le \, p$.
Rewriting it in the form $\sin t\,  \le \, p(1-\cos t) = 2p\sin^2(t/2)$, substituting~$\sin t = 2\sin (t/2) \cos (t/2)$
and dividing both sides by~$2\sin (t/2)$ we come, to an equivalent inequality~$\cos (t/2) \le p \sin (t/2)$,
or $\tan (t/2) \ge \frac{1}{p}$. For~$t \in \bigl( 0, \frac{\pi}{2}\bigr]$, this inequality is equivalent to
$$
\sin t \ = \  \frac{2 \tan (t/2)}{1 + \tan^2 (t/2)} \ \ge \ \frac{2 p}{1 + p^2}\, ,
$$
 which follows from the
assumption~$\sin t \ge \frac{2}{p}$.

$\mathbf{2).\  \sin t  < \frac{2}{p}}.$ In this case we also establish a stronger inequality:
$\sin t + p \le \sqrt[3]{p^3 + \frac{20}{\sin t}}$. Taking the third power of the
both sides and writing~$\sin t = s$, we get after simplifications
$$
s^4 \, + \, 3 p s^3\, + \, 3 p^2 s^2 \ \le \ 20\, .
$$
The left hand side is increasing in~$s$, therefore, to establish this inequality
for all~$s \in \bigl(0,  \frac{2}{p}\bigr)$ it suffices to check it at $s = \frac{2}{p}$, where
 it becomes~$\frac{16}{p^4} + \frac{24}{p^2} + 12 \, \le \, 20$. This holds for all~$p \ge 2$.

{\hfill $\Box$}
\medskip

{\tt Proof of Theorem~\ref{th20}.}
{\bf 1) Bound on the rate of growth.}
The matrices in (\ref{pr}) have a block upper triangular form with two diagonal blocks of dimensions~$1$ and~$2$. Let us denote
\begin{equation}\label{pr1}
P \ = \
\left(
\begin{array}{ll}
1 & 0 \\
0 & 0
\end{array}
\right)\, , \quad
R \ = \
\left(
\begin{array}{ll}
\cos \alpha & -\sin \alpha\\
\sin \alpha & \ \, \cos \alpha
\end{array}
\right)\, , \quad a \ = \ (\sin \alpha \ , \ \cos \alpha - 1)^T\, .
\end{equation}
Thus, $A_0^{(2, 2)} = P$ is an orthogonal projection onto the first coordinate axis, $A_1^{(2, 2)} = R$
is a rotation of angle $\alpha$ around the origin.
We have~$A_1^{(1, 2)} = (\sin \alpha \, , \cos \alpha - 1)$, denote this $1\times 2$-matrix by~$Q$.
Thus, $Q = a^T$.

For every~$n \in \n$,  we have~$A_0^n = A_0$,  $(A_1^n)^{(1,1)} = 1$, and
$\, (A_1^n)^{(2,2)} = R^n$.
We now consider the block~$(A_1^n)^{(1,2)}$ which will be denoted by~$Q_n$.
 We also write~$e_1, e_2$ for~the basis vectors of~$\re^2$.
For our purposes,  it suffices to evaluate~$Q_ne_1$, which is the scalar product of the vectors~$Q_n^T$ and~$e_1$. For a vector~$z = (x, y)^T$, we denote by~$\mathbf{z} = x + i y$ the corresponding complex number.
We have
$$
\langle a , z\rangle \, = \, x \sin \alpha \, + \, y (\cos \alpha - 1) \ = \
{\rm Re} \, \bigl(\, - i \, (e^{i \alpha} - 1) \, \mathbf{z}\bigr).
$$
The complex number associated with the vector~$R^n e_1$ is~$e^{i n \alpha}$, in particular, $\mathbf{e_1} = 1$. Therefore,
$$
Q_n e_1 \ = \ \sum_{l=0}^{n-1} \langle a\, , \, R^l e_1 \rangle \ = \
\,  \sum_{l=0}^{n-1} {\rm Re}\, \left( - i \, (e^{i \alpha} - 1)\, e^{i \, l \alpha}\, \right) \ = \
$$
$$
\ = \ {\rm Re}\, \left(  - i \, (e^{i \alpha} - 1)\, \frac{e^{i n \alpha} - 1}{e^{i \alpha} - 1}\right)\ =
\  {\rm Re}\, \Bigl( - i \, (e^{i n \alpha} - 1)\Bigr) \ = \ \sin n \alpha\, .
$$
Thus, $Q_n e_1 = \sin n \alpha$. Finally, $PR^ne_1$ is an orthogonal
projection of the vector~$R^ne_1$ on the first coordinate axis. Whence, $PR^ne_1 \, =\, (\cos n\alpha )\, e_1$.

Now we estimate the norm of an arbitrary product~$\Pi \in \cA^N$. First of all, since the matrix~$A_1$ has
three simple eigenvalues of modulus~$1$, it follows that~$\|A^N\| \le C$, where $C$ does not depend on~$N$.
Whence, if the product~$\Pi$ does not contain~$A_0$, then its norm is bounded by~$C$ uniformly for all~$N$.
Moreover, if~$\Pi$ ends with some power of~$A_1$,
one can remove this power from the product, and this may reduce the norm of~$\Pi$ by at most~$C$ times
(because of the submultiplicativity of the norm).  In view of this, since our goal is to bound the norm of an arbitrary product~$\Pi,$ we can restrict our attention to such products ending with~$A_0$. We have
$$
\Pi \ = \ A_1^{n_1}A_0 A_1^{n_{2}}A_0 \cdots A_1^{n_{k-1}}A_0A_1^{n_k}A_0\, .
$$
Since~$\Pi^{(1,1)} = 1$ and~$\|\Pi^{(2,2)}\| \le 1$, we only need to estimate~$\|\Pi^{(1,2)}\|$.
Moreover, $\Pi$ ends with~$A_0$, hence $\|\Pi^{(1,2)}\| = |\Pi^{(1,2)}e_1|$.
Since~$A_0^{(1,2)} = 0$, it follows that
$$
\Pi^{(1,2)}\, e_1 \ = \ Q_{n_k}e_1 \, + \, Q_{n_{k-1}}P R^{n_k}e_1 \, + \, \cdots \,  + \, Q_{n_1}P R^{n_{2}}P \cdots R^{n_{k-1}}PR^{n_k}e_1.
$$
Applying the above equalities~$PR^n e_1 = (\cos n\alpha ) \, e_1\, $ and~$\, Q_ne_1 = \sin n \alpha$, we obtain
\begin{equation}\label{main}
\Pi^{(1,2)}\, e_1 \ = \  \sin n_k \alpha \ + \ \sin n_{k-1} \alpha \ \cos n_k \alpha \ + \
\cdots \ + \
 \sin n_{1}\alpha \ \cos n_{2}\alpha \ \cdots \ \cos n_{k}\alpha.
\end{equation}
Now, replacing all these sines and cosines by their
modules does not reduce the modulus of this sum. Consider a sequence of numbers~$\{S_r\}_{r=1}^{k}$
defined by the following recursion:
$$
S_0 = 0\, , \qquad  S_{r} \ = \ |\sin n_r \alpha | \  + \ S_{r-1}|\cos n_{r} \alpha |\, , \quad r = 1, \ldots , k\, .
$$
If we removed modules in this recurrence relation, we would get sum~(\ref{main}) after the $k$th iteration.
Therefore, $S_k$ is greater than or equal to the modulus of this sum.
Thus, $S_k \ge |\Pi^{(1,2)}\, e_1|$.  To estimate~$S_k$ we now invoke  Lemma~\ref{l30}.
 Denote by~$r_0$ the greatest number~$r \in \{0, \ldots , k\}$ for which~$S_r < 8$.
 If~$r_0 =k$, then~$S_k < 8$. Otherwise, if $r_0 \le k-1$, then
  $S_{r_0+1} \le 1 + S_{r_0} < 9$. In particular, if $r_0 = k-1$, then~$S_k < 9$.
  Thus, without loss of generality it can be assumed that~$r_0 \le k-2$.
Applying Lemma~\ref{l30} consecutively for~$r = r_0+2, \ldots , k$, substituting~$p = \sqrt[3]{S_{r-1}}$
and~$t \in \bigl(0, \frac{\pi}{2}\bigr]$
such that $\sin t = |\sin n_{r}\alpha|\, , \, \cos t = |\cos n_{r}\alpha|$, we obtain
\begin{equation}\label{itog1} S_k \ \le \ \left( S_{r_0+1} \ + \ \sum_{r=r_0+2}^{k} \, \frac{20}{|\sin n_r \, \alpha|}\, \right)^{\, 1/3}  \ < \quad
  \left( 9 \ + \ \sum_{r=1}^{k} \ \frac{20}{|\sin n_r \, \alpha |}\, \right)^{\, 1/3}  \, .
\end{equation}
Now it is time to invoke our assumption that $\alpha/\pi$ is a quadratic irrational.
Denote by~$h(x)$ the distance from a number~$x \in \re$ to the closest integer.
For every~$t\in \re$, we have~$|e^{2i t} - 1| \ge 4 h \bigl(\frac{t}{\pi}\bigr)$.
By the Liouville approximation theorem (see~\cite{apostol})
there exists a constant $m_{\alpha}$ such that $\bigr|\,  \frac{\alpha }{\pi} - \frac{p}{q}\, \bigr| \ge \frac{m_{\alpha}}{q^2}$, for all integers~$q, p$. Consequently,
$h\, \bigl(n\frac{\alpha }{\pi}\bigr) \ge \frac{m_{\alpha}}{n}\, , \ n \in \n$.
Therefore, for each~$r$, we have
$$
   \frac{1}{|\sin n_r \, \alpha|} \ = \ \frac{1}{|e^{2 i n_r \, \alpha } - 1|} \ \le \
   \frac{1}{4 h\, \bigl(\frac{n_r \alpha }{\pi}\bigr)}\ \le \ \frac{n_r}{4 m_{\alpha}}\, .
$$
Substituting to~(\ref{itog1}), we obtain
$$
\, S_k \, \le \, \left(\, 9 \, + \, \frac{5}{ m_{\alpha}}\, \sum_{r=1}^k n_r \, \right)^{\, 1/3} \ \le \
\left(\, 9 \, + \, \frac{5}{ n_{\alpha}}N \right)^{\, 1/3}\ \le \ C_{\alpha}\, N^{1/3},
$$
where the constant~$C_{\alpha}$ does not depend on~$N$. This completes the proof of the upper bound~(\ref{cubic}).
\smallskip

{\bf 2) Existence of an unbounded product.}
Now let us  show that this upper bound is sharp. First, we prove that there are arbitrarily
long products~$\Pi \in \cA^N$ such that~$\|\Pi\| \ge C_0N^{1/3}$, where $C_0$ is a constant depending only on~$\alpha$.
Then we use those products to build an infinite product of matrices from~$\cA$ for which there exists
a sequence~$\{N_j\}_{j\in \n}$ such that $\lim\limits_{j \to \infty}\, \frac{\log \|\Pi_j\|}{\log N_j}\, \to \, \frac{1}{3}$,
where $\Pi_j$ is the suffix of length~$N_j$ of the infinite product.

We start by constructing finite products~$\Pi \in \cA^N$ such that $\|\Pi\| \ge C_0 N^{1/3}$.
To this end we take $N = n^3+n^2, k = n^2, \, n_1 = \ldots = n_k = n$, where the number~$n \in \n$
is such that
\begin{equation}\label{n}
\left\{\, n \, \frac{\alpha}{2\pi}\right\} \ \le \ \frac{1}{n}\, .
\end{equation}
(where~$\{x\} = x - [x]$ represents the fractional
part of~$x$). By Kronecker's approximation theorem \cite[Chapter 7]{apostol} , there are infinitely many such natural numbers~$n$, provided
$\frac{\alpha}{2\pi}$ is irrational. For the sake of simplicity, we assume~$\alpha > 0$. We take
$\Pi = (A_0A_1^{n})^{n^2}$ and show that $\|\Pi\| \ge C_0 N^{1/3}$.
Equality~(\ref{main}) now reads
\begin{equation}\label{cubic2}
\Pi^{(1,2)}\, e_1 \ = \  \sin n \alpha \, \left( \sum_{r=0}^{n^2-1} (\cos\, n\alpha)^{r}\right)\ = \
\sin n \alpha \ \frac{ 1- (\cos\, n\alpha)^{n^2}}{1 - \cos \, n\alpha }\, .
\end{equation}
Since $\frac{m_{\alpha}}{2n}\le \bigl\{n \frac{\alpha}{2\pi}\bigr\} \le \frac{1}{n}$ it follows that
there exists $t \in \bigl(\frac{m_{\alpha}\pi}{n}\, , \, \frac{2\pi}{n}\bigr)$ such that  $t  \equiv n \alpha\pmod {2\pi}$.
Therefore, $\ \sin n \alpha \, \sim \, t\, $ and $ \, 1 - \cos n \alpha \, \sim \,\frac{t^2}{2}$ as~$n \to \infty$.
Furthermore,
\begin{equation}\label{cubic3}
1 \ - \ (\cos\, n\alpha)^{\, n^2} \ \sim \ \frac{n^2t^2}{2} \ \ge \ \frac{m_{\alpha}^2\pi^2}{2}\, .
\end{equation}
Substituting in~(\ref{cubic2}), we get for large~$n$:
$$
\Pi^{(1,2)}\, e_1 \ \ge \ \frac{t\, m_{\alpha}^2\pi^2/2}{t^2/2} \ = \ \frac{\, m_{\alpha}^2\pi^2}{t} \ \ge \
n \, \frac{\, \, m_{\alpha}^2\pi^2}{2\pi} \ = \ n \, \frac{ \, m_{\alpha}^2 \pi}{2} \, .
$$
Note that $N^{1/3} = (n^3 + n^2)^{1/3} \le 2n$ for all~$n \in \n$. Whence, $n \ge \frac{1}{2}N^{1/3}$, and
we finally obtain
\begin{equation}\label{cubic4}
\Pi^{(1,2)}\, e_1 \ \ge  \  \frac{\, m_{\alpha}^2\pi}{4}\, N^{1/3} \, .
\end{equation}
Thus, for every product of the form~$\Pi = (A_0A_1^{n})^{n^2}$, where $n$ is an arbitrary number satisfying~(\ref{n}), we have $\|\Pi\| \, \ge \, C_0\, N^{1/3}$, where $\ C_0 = \frac{\, m_{\alpha}^2\pi}{4}$.
\smallskip

Now we build an infinite product~$\Pi$ with blocks of this form that has the exponent of growth~$1/3$.
We take a sequence of numbers~$\{n_j\}_{j \in \n}$ that satisfy~(\ref{n}). The growth of this sequence can be
chosen arbitrarily fast and will be specified later. To each~$n_j$ we assign the corresponding
product~$P_j = (A_0A_1^{n_j})^{n_j^2}$
of length~$N_j = n_j^3 + n_j^2$ and a number~$t_j = n_j \alpha$ such that
$\, t_j \, \equiv \, n_j \alpha\, \pmod {2\pi}$ and
$t_j \in \bigl(\frac{m_{\alpha}\pi}{n_j}\, , \, \frac{2\pi}{n_j}\bigr)$.

Denote~$\Pi_k = \prod_{j=1}^k P_j$. For every~$j$, we have
$P_j^{(1,2)}e_1 \ge C_0 N_j^{1/3}\, $ and $\, P_j^{(2,2)} = (\cos\, t_j )^{n_j^2} \, \ge \, C_1$,
where $C_1 > 0$ does not depend on~$j$. Indeed, since $\ln (1-x) > -2x$ for $x \in \bigl(0, \frac12 \bigr)$, we have
$$
(\cos\, t_j )^{n_j^2} \ \ge \ \left(1 - \frac{t_j^2}{2} \right)^{n_j^2} \ = \ e^{n_j^2\ln \bigl(1 - \frac{t_j^2}{2} \bigr)}\ \ge \ e^{-n_j^2 t_j^2} \ \ge \ e^{-4\pi^2}\, .
$$
It remains to set~$C_1 = e^{-4\pi^2}$. Therefore, $\Pi_k e_1 \, \ge \, C_0 \, \Bigl(\, \sum_{j=1}^k C_1^{j-1}N_j^{1/3} \, \Bigr)$.
If $N_1, \ldots , N_{k-1}$ are already fixed, we take the next length~$N_k$ so that
$$
C_0 \left(\sum_{j=1}^k C_1^{j-1}N_j^{\, \frac{1}{3}} \right)\ \ge \ \left(\sum_{j=1}^k N_j\, \right)^{\frac{1}{3} - 2^{-k}} \, .
$$
This is possible, because the left hand side grows as $N_j \to \infty$ faster than the right hand side.
Writing~$M_k$ for the length of~$\Pi_k$, we obtain $\frac{\log \|\Pi_k\|}{\log M_k}\, \to \, \frac13$ as $k \to \infty$, which concludes the proof.

\section{Open problems}

In order to end this paper, we propose here a series of natural open questions coming out of our work.  The first one can roughly be stated as ``What is the minimal worst case rate of growth of a marginally stable semigroup of matrices? ''  We can make this question formal in two different ways:

\begin{opq}
\begin{itemize}
\item (Weak version) Does there exist a marginally unstable set of matrices with growth slower than $x^\alpha$ for any $\alpha>0$ ? That is, the set is unstable, but for any trajectory $x(\cdot),$
$$\lim_{t \rightarrow \infty}{\log{|x(t)|}/\log{t}}=0? $$
\item (Strong version) Does there exist a marginally unstable set of matrices with arbitrarily slow growth? That is, for any increasing sequence of numbers $a_1<a_2<\ldots$ such that $a_t \to \infty\, , \, t \to \infty$ there exists an unstable set for which every trajectory $x(\cdot)$ satisfies
$$x(t)<a_t.$$
\end{itemize}
\end{opq}

Another interesting open question, in our view, is about matrices leaving a common cone invariant.
In recent years, much research effort has been devoted to study sets of matrices that leave a certain cone invariant. Indeed, it has been progressively realized that these sets enjoy much stronger properties than arbitrary sets of matrices.  Often, these properties allow one to design efficient algorithmic procedures in order to answer questions which are known to be hard (or even impossible) to answer in general.
Examples of such problems include stability or stabilizability \cite{conic}, continuity of joint spectral characteristics\cite{jungers_invariant},  see also \cite{valcher_misra,fainshil2009stability}.

For these reasons, it seems interesting to study sets of matrices sharing an invariant cone, or more specifically, the case where this invariant cone is the positive orthant.
\begin{opq}
What is the minimal worst case rate of growth of marginally unstable positive systems?  (i.e. characterized by nonnegative matrices in discrete time, Metzler matrices in continuous time.)
\end{opq}
Our interest in this particular subcase is twofold: it might be easier to answer our open questions in this case, or a lower bound on the minimal growth might be larger than for general systems.

Finally, as the reader can see, our Section \ref{sec-sublin} is exclusively devoted to discrete time systems.  This leads us to the obvious question:

\begin{opq}
Are there switching systems in continuous time which are marginally unstable, but with sublinear growth?
\end{opq}

{\hfill $\Box$}
\section*{Acknowledgements}
The research was carried out
when the first  author was visiting the department of Mathematical Engineering,
Universit\'e Catholique de Louvain (UCL), Belgium. He is very grateful to the university for hospitality.

\end{document}